\documentclass[12pt]{amsart}
\usepackage{amsmath,amssymb,amsthm}

\usepackage{color}

\newcommand{\R}{\mathbb R}

\newcommand{\E}{\mathbb E}

\newtheorem{thm}{Theorem}[section]
\newtheorem{cor}[thm]{Corollary}
\newtheorem{lem}[thm]{Lemma}
\newtheorem{prop}[thm]{Proposition}
\theoremstyle{definition}
\newtheorem{defn}[thm]{Definition}
\theoremstyle{remark}
\newtheorem{rem}[thm]{Remark}
\newcommand{\ds}{\displaystyle}
\begin{document}

\title[On the Invariant Theory of Weingarten Surfaces]
{ON THE INVARIANT THEORY OF WEINGARTEN SURFACES IN EUCLIDEAN SPACE}%
\author{Georgi Ganchev and Vesselka Mihova}%
\address{Bulgarian Academy of Sciences, Institute of Mathematics and Informatics,
Acad. G. Bonchev Str. bl. 8, 1113 Sofia, Bulgaria}%
\email{ganchev@math.bas.bg}%
\address{Faculty of Mathematics and Informatics, University of Sofia,
J. Bouchier Str. 5, 1164 Sofia, Bulgaria}
\email{mihova@fmi.uni-sofia.bg}
\subjclass[2000]{Primary 53A05, Secondary 53A10}%
\keywords{Strongly regular surfaces, Weingarten surfaces, geometric principal
parameters, canonical principal parameters}%

\begin{abstract}
We prove that any strongly regular Weingarten surface in Euclidean space
carries locally geometric principal parameters. The basic theorem states
that any strongly regular Weingarten surface is determined up to a motion
by its structural functions and the normal curvature function satisfying
a geometric differential equation. We apply these results to the special
Weingarten surfaces: minimal surfaces, surfaces of constant mean
curvature and surfaces of constant Gauss curvature.
\end{abstract}
\maketitle
\section{Introduction}

Our aim in this paper is to introduce and prove a local existence of
geometric principal parameters on a class of surfaces in Euclidean space.

For a surface $S$ in the Euclidean space ${\R}^3$ we use the notations
$E,F,G; \; L,M,N$ for the coefficients of the first and the second
fundamental form, respectively. $K$ and $H$ will stand for the Gauss
and mean curvature, respectively.

A useful tool for investigating surfaces in Euclidean space are the
isothermal parameters, which are characterized by the conditions
$E=G, \; F=0$.

In \cite{R} Rellich considered isothermal parameters $(u,v)$ on a
surface $S$ and proved that $S$ is uniquely determined by the
(natural) equation $H=f(u,v)$. Similar questions are studied
with respect to the second fundamental form and the (natural)
equation $K=f(u,v)$. In the latter case the appropriate parameters
satisfy the condition $L=N \neq 0, \; M=0$ and the
equation $K=f(u,v)$ determines uniquely the surface $S$.

In \cite{S1} Scherrer studied the question of determining a surface
$S$ with:

(i) the coefficients of the first fundamental form I and the mean
curvature $H$;

(ii) the coefficients of the second fundamental form II and the
Gauss curvature $K$;

(iii) the coefficients of the third fundamental form III and the
harmonic curvature $\ds{\frac{K}{H}}$.

In any of these cases it is proved a theorem of existence and uniqueness
of a surface $S$.

Looking for a natural equation of a surface
$S: x=x(u,v), \; (u,v)\in \mathcal D$,
Scherrer introduced in \cite{S2} the functions

(i) \emph{support} function: $p=lx$, where $l$ is the unit
normal field to $S$;

(ii) \emph{radius} function: $r=\sqrt{x^2}$

\noindent
and proved that the functions $p, \, r$ and any of the functions
$H,\,K$ or $\ds{\frac{K}{H}}$ determine $S$ up to a motion.

In \cite{L} Leichtweiss considered \emph{normalized} isothermal
parameters $(p,q)$ on a surface $S$ characterized by the conditions:

(i) any of the fundamental forms I,\,II or III can be written in the form
$D(dp^2+dq^2)$;

(ii) an arbitrary given analytic curve on $S$ (with natural parameter $s$)
is determined by the equalities ($p=s,\;q=0$).

Then in a stripe-neighborhood of the point ($p=0,\,q=0$)
(determined by the given curve) the surface $S$ is uniquely
determined by a linear-fractional function of $H$ and $K$.

In \cite{B} Bryant introduces natural principal coordinates, which are suitable for
integrating the system of Gauss-Codazzi equations of a surface with prescribed shape
operator.

In this paper we consider surfaces in Euclidean space with respect to principal
parameters. For any surface $S: \; x=x(u,v), \; (u,v)\in \mathcal D$
we denote by $\nu_1,\, \nu_2; \; \gamma_1, \, \gamma_2$ the principal normal
curvatures and the principal geodesic curvatures, respectively.

We define the class of \emph{strongly regular} surfaces by the condition
$$[\nu_1(u,v)-\nu_2(u,v)] \gamma_1(u,v) \gamma_2(u,v)\neq 0, \quad (u,v)\in \mathcal D.$$

In Theorem \ref{T:2.1} we reformulate the fundamental Bonnet theorem for strongly regular
surfaces in terms of the four invariants $\nu_1,\, \nu_2; \; \gamma_1, \, \gamma_2$.

For the determining conditions of the Weingarten surfaces we use the
following two forms:
$$\nu_2=w(\nu_1) $$
or
$$\nu_1=f(\nu), \quad \nu_2=g(\nu).$$

In Theorem \ref {T:4.2} we prove that any strongly regular Weingarten surface
admits locally geometric principal parameters.

In the main Theorem \ref {T:4.3} we prove that:

{\it Any strongly regular Weingarten surface is determined up to a motion
by its determining functions $f$, $g$ and the function
$\nu$, satisfying the geometric partial differential equation.}

As a consequence we obtain the corresponding theorems for the classes
of surfaces determined by the conditions: $H=0, \; H=1/2, \; K=1, \; K=-1$.

\section{Surfaces one of whose family of principal lines consists
of geodesics}

Let $S: \, x=x(u,v), \; (u,v) \in {\mathcal D}$ be a surface in
${\R}^3$ with the standard metric and its flat Levi-Civita connection
$\nabla$. The unit normal to $S$ is denoted by $l$ and
$E, F, G; \; L, M, N$ stand for the coefficients of the first and
the second fundamental forms, respectively. We suppose that $S$ has
no umbilical points and the families $\mathcal F_1$ and
$\mathcal F_2$ of the principal lines on $S$ form a parametric net, i. e.
$$F(u,v)=M(u,v)=0, \quad (u,v) \in \mathcal D.$$
Then the principal curvatures $\nu_1$ and $\nu_2$ are
$$\nu_1=\frac{L}{E}, \qquad \nu_2=\frac{N}{G}.$$
The principal geodesic curvatures (geodesic curvatures of the principal
lines) $\gamma_1$ and $\gamma_2$ of $S$ are
$$\gamma_1=-\frac{E_v}{2E\sqrt G}, \qquad \gamma_2=
\frac{G_u}{2G\sqrt E}, \leqno(2.1)$$
respectively.

We define the tangential frame field $\{X, Y\}$ as follows
$$X:=\frac{x_u}{\sqrt E}, \qquad Y:=\frac{x_v}{\sqrt G}$$
and assume that the moving frame field $XYl$ is always right oriented orthonormal.
The following Frenet type formulas for $XYl$ are valid

$$\begin{tabular}{ll}
$\begin{array}{llccc}
\nabla_{X} \,X & = &  &\gamma_1 \,Y + \nu_1 \, l,  &\\
[2mm]
\nabla_{X} Y & = -\gamma_1 \, X, & & \\
[2mm]
\nabla_{X} \, l & = - \nu_1 \, X, & & &
\end{array}$ &
\quad
$\begin{array}{llccc}
\nabla_{Y} \,X & = & & \gamma_2 \, Y, &\\
[2mm]
\nabla_{Y} Y & = -\gamma_2 \, X & &  & + \nu_2 \, l,\\
[2mm]
\nabla_{Y}\, l & = &  & - \nu_2 \, Y. &
\end{array}$
\end{tabular}\leqno(2.2)$$
\vskip 2mm
Then the Lee bracket $[X,Y]$ is given by $[X,Y]=-\gamma_1\,X-\gamma_2\,Y$.

The integrability condition
$\nabla_X\nabla_Yl-\nabla_Y\nabla_Xl-\nabla_{[X,Y]}l=0$
is equivalent to the Codazzi equations
$$\begin{array}{l}
\displaystyle{\gamma_1=\frac{Y(\nu_1)}{\nu_1-\nu_2}=
\frac{(\nu_1)_v}{\sqrt G(\nu_1-\nu_2)}},\\
[4mm]
\displaystyle{\gamma_2=\frac{X(\nu_2)}{\nu_1-\nu_2}=
\frac{(\nu_2)_u}{\sqrt E(\nu_1-\nu_2)}}
\end{array}\leqno(2.3)$$
and the integrability condition
$\nabla_X\nabla_YY-\nabla_Y\nabla_XY-\nabla_{[X,Y]}Y=0$
implies the Gauss equation
$$Y(\gamma_1)-X(\gamma_2)-(\gamma_1^2+\gamma_2^2)=\nu_1\nu_2=K,$$
or
$$\frac{(\gamma_1)_v}{\sqrt G}-\frac{(\gamma_2)_u}{\sqrt E}-
(\gamma_1^2+\gamma_2^2)=\nu_1\nu_2=K.\leqno(2.4)$$

Formulas (2.2) imply explicit expressions for the curvature and
the torsion of any principal line.

Let $c_1: z=z(s), \; s \in J$ be a line from the family $\mathcal F_1$
($v={\rm const}$) parameterized by a natural parameter and
$\kappa_1, \, \tau_1$ be its curvature and torsion,
respectively. Since $c_1$ is an integral line of the unit vector
field $X$, then
$$z'=X, \quad z''=\nabla_XX=\gamma_1\,Y+\nu_1\,l, \quad
z'''=\nabla_X\nabla_XX=-(\gamma_1^2+\nu_1^2)\,X+X(\gamma_1)\,
Y+X(\nu_1)\,l,$$
$$\kappa_1^2= \gamma_1^2+\nu_1^2.\leqno(2.5)$$
Assuming that $\gamma_1^2+\nu_1^2>0$ along $c_1$, we find
$$\tau_1=\frac{\gamma_1\, X(\nu_1)-\nu_1\,X(\gamma_1)}{\gamma_1^2+
\nu_1^2}= -X\left(\arctan\frac{\gamma_1}{\nu_1}\right)=
-X(\theta_1),\leqno(2.6)$$
where $\theta_1=\angle (n_1, l)$ is the angle between the principal
normal $n_1$ of $c_1$ and the normal $l$ of the surface  $S$.

For the lines $c_2$ of the family $\mathcal F_2$ we obtain in a
similar way the corresponding formulas
$$\kappa_2^2=\gamma_2^2+\nu_2^2,\leqno(2.7)$$
$$\tau_2=\frac{\gamma_2\, Y(\nu_2)-\nu_2\,Y(\gamma_2)}
{\gamma_2^2+\nu_2^2}=-Y\left(\arctan\frac{\gamma_2}{\nu_2}\right)=
-Y(\theta_2), \quad  (\gamma_2^2+\nu_2^2>0).\leqno(2.8)$$
\vskip 2mm
The aim of this section is to describe locally in a geometric
(constructive) way the class of surfaces whose family
$\mathcal F_1$ consists of geodesics.

Let $c_2: x=x(v), \; v\in J_2,$ be a smooth regular curve in
${\E}^3$ parameterized by a natural parameter $v$ with vector
invariants $t(v), \, n(v), \, b(v)$,  curvature
$\kappa (v)>0$ and torsion $\tau (v)$.

A unit normal vector field $e(v)$ along the curve $c_2$ is said to be
\emph{torse-forming} \cite{Y} if $\nabla_te=ft$ \,for a certain function
$f(v), \; v\in J_2$ on the curve $c_2$.

We shall consider an orthonormal pair $\{e_1(v), e_2(v)\}$ of
torse-forming normals along the curve $c_2$. Denoting by
$\theta=\angle(n,e_1)$, then the pair
$$\begin{array}{l}
e_1=\cos \theta \, n +\sin \theta \, b,\\
[2mm]
e_2=-\sin \theta \, n + \cos \theta \, b,\end{array}$$
is determined uniquely up to a constant angle $\theta_0$ by the condition
$\displaystyle{\theta(v)=-\int_0^v\tau\,dv+\theta_0}$.

Further we choose $\theta_0=0$, i.e. the pair $\{e_1(v), e_2(v)\}$
satisfies the initial conditions $e_1(0)=n(0), \, e_2(0)=b(0)$.

The orthonormal frame field $t(v)e_1(v)e_2(v)$ satisfies the Frenet
type formulas
$$\begin{array}{l}
t'=\qquad \qquad \quad \kappa \, \cos \theta\, e_1 - \kappa \,
\sin \theta \, e_2,\\
[2mm]
e_1'=- \kappa \, \cos \theta \, t,\\
[2mm]
e_2'= \; \; \; \kappa \, \sin \theta \, t.\end{array}\leqno(2.9)$$

We consider the regular plane curve
$$c_1: z(u)=x(v) + \lambda(u)\, e_1(v)+ \mu(u)\, e_2(v) \quad u
\in J_1,$$
rigidly connected with every Cartesian coordinate system
$x(v)e_1(v)e_2(v)$. We suppose that the parameter $u$ is natural
($\dot \lambda^2 + \dot \mu^2=1$) and the functions
$\lambda, \, \mu$ satisfy the initial conditions
$$\begin{array}{l}
\lambda(0)=\mu(0)=0;\\
[2mm]
\dot\lambda(0)=1, \; \dot \mu(0)=0.\end{array}$$
Then the (plane) curvature $\kappa_1=\kappa_1(u)\neq 0$ of $c_1$
completely determines the functions $\lambda$ and $\mu$.
\vskip 2mm
Now, let us consider the surface
$$S: \; Z(u,v)= x(v) + \lambda(u)\, e_1(v)+\mu(u)\, e_2(v);
\quad u \in J_1, \; v \in J_2.\leqno(2.10)$$
\vskip 2mm
Taking into account (2.9) and (2.10) we obtain
$$\begin{array}{l}
Z_u=\dot \lambda \, e_1+\dot \mu \, e_2,\\
[2mm]
Z_v=[1-\kappa(\lambda \cos \theta - \mu \sin \theta)]\,t,\\
[2mm]
Z_u\times Z_v=-[1-\kappa(\lambda \cos \theta - \mu \sin \theta)]
(-\dot \mu \, e_1+ \dot \lambda \, e_2).\end{array}$$

The surface $S$ is smooth at the points, where
$$\lambda \cos \theta - \mu \sin \theta \neq \frac{1}{\kappa}.
\leqno(2.11)$$
We orientate the surface $S$ by choosing $l=-\dot \mu \, e_1+
\dot \lambda \, e_2$, i. e. the normal to $S$ is the plane normal to
$c_1$. Then
$$\begin{array}{l}
E=1, \quad F=0, \quad G=[1-\kappa(\lambda \cos \theta -
\mu \sin \theta)]^2;\\
[2mm] L=\kappa_1, \quad M=0, \quad N= -\kappa\,(\dot\lambda\sin
\theta+\dot\mu\cos\theta) [1-\kappa(\lambda \cos \theta - \mu \sin
\theta)].\end{array}\leqno(2.12)$$

We denote by $\Gamma$ the class of surfaces, given by $(2.10)$,
under the condition $(2.11)$.

Taking into account (2.1) and (2.12), we obtain:
\begin{lem}
Any surface $S$ of the class $\Gamma$ has the following properties:

$1)$ the parametric lines are principal;

$2)$ the family $\mathcal F_1$ consists of geodesics.
\end{lem}
The invariants of any surface from the class $\Gamma$ are given as follows:
$$\begin{array}{ll}
\nu_1=\kappa_1(u), & \ds{\nu_2=\frac{-\kappa\,(\dot\lambda\sin \theta+
\dot\mu\cos\theta)}
{1-\kappa(\lambda \cos \theta - \mu \sin \theta)};}\\
[4mm]
\gamma_1= 0, & \ds{\gamma_2=\frac{-\kappa (\dot \lambda \cos \theta -
\dot \mu \sin \theta)}
{1-\kappa(\lambda \cos \theta - \mu \sin \theta)}\,.}\end{array}
\leqno(2.13)$$
\begin{thm}\label{T:1.1}
Let $S$ be a surface parameterized by principal parameters.
If the family $\mathcal F_1$ of principal lines consists of regular
geodesics, then $S$ is locally part of a surface from the class $\Gamma$.
\end{thm}
\emph{Proof:} Let us consider the derivative formulas
$$\begin{array}{l}
\nabla_XX= \quad \quad \quad \gamma_1\,Y+\nu_1 \,l,\\
[2mm]
\nabla_XY=-\gamma_1\,X,\\
[2mm]
\nabla_X\,l\;=-\nu_1\,X.
\end{array}\leqno(2.14)$$

Under the condition $\gamma_1=0$ it follows from (2.14), (2.3) and
(2.6) that:

\begin{itemize}
  \item [-]{any line $c_1(v)$ from the family $\Gamma_1$ lies in a plane,
  which we denote by $E^2(v)$;}
  \item [-]{the function $\nu_1(u,v)=\nu_1(u)$ depends only on $u$ and
  $\kappa_1=\nu_1(u)\neq 0$;}
  \item [-]{the normal $l$ of $S$ lies in $E^2(v)$, i. e. $l$ is
  (up to a sign) the plane normal of $c_1$.}
\end{itemize}

Choose one of the principal lines from $\mathcal F_2$ as the line
$c_2: \, x=x(v)$. According to (2.2), the unit vector fields $X$ and
$l$ are torse-forming normals along $c_2$. Considering the parametric
line $c_1$ in the plane $x(v)Y(v)l(v)$, $v$ - fixed, we obtain that
$c_1$ is uniquely determined by its curvature function $\kappa_1=\nu_1$
and the initial Frenet frame $x(v)Y(v)l(v)$. Denoting
$e_1(v)=Y(v), \, e_2(v)=l(v)$, we obtain that the curve $c_1$ is rigidly
connected with the frame field $x(v)e_1(v)e_2(v)$. We conclude from here
that the surface $S$ has locally the construction of a surface from the
class $\Gamma$. \qed

\section{On the invariant theory of strongly regular surfaces}

In this section we introduce the class of strongly regular surfaces in
Euclidean space. For this class we reformulate the fundamental
Bonnet theorem in terms of the invariants of the surfaces.

\begin{defn}
A surface $S: \; x=x(u,v), \; (u,v)\in \mathcal D$ is said to be \emph{strongly regular}
if
$$[\nu_1(u,v)-\nu_2(u,v)]\gamma_1(u,v)\gamma_2(u,v) \neq 0, \quad (u,v)\in\mathcal D.$$
\end{defn}

The following statement follows in a straightforward way.

\begin{lem} If $S$ is a strongly regular surface, then the sign of the function
$(\nu_1-\nu_2)\gamma_1\gamma_2$ does not depend on the parametrization of $S$.
\end{lem}

In what follows we specialize the right oriented moving frame $XYl$ so that $\nu_1-\nu_2>0$.
Then the class of strongly regular surfaces is divided into two subclasses
characterized with the following conditions
$$ \gamma_1 \gamma_2 > 0, \qquad \gamma_1 \gamma_2 < 0,$$
respectively.

Considering Codazzi formulas (2.3) we make the observation that
$$\begin{array}{l}
\gamma_1 \neq 0 \quad \iff \quad (\nu_1)_v\neq 0;\\
[2mm]
\gamma_2 \neq 0 \quad \iff \quad (\nu_2)_u\neq 0.\end{array}$$

Then it follows that for any strongly regular surface
$$\sqrt E=\frac{(\nu_2)_u}{\gamma_2(\nu_1-\nu_2)} >0, \quad
\sqrt G=\frac{(\nu_1)_v}{\gamma_1 (\nu_1-\nu_2)}>0\, .\leqno(3.1)$$

Thus, using (2.3) we can express the coefficients
$E,\,G$ and $L,\,N$ by the invariant functions $\nu_1,\,\nu_2$ and
$\gamma_1,\,\gamma_2$. In view of (3.1) the Frenet formulas (2.2) get
the form
$$\begin{array}{lccc}
X_u  = &  & \displaystyle{\frac{\gamma_1 \, (\nu_2)_u}
{\gamma_2(\nu_1-\nu_2)}\, Y}
&+ \, \displaystyle{\frac{\nu_1 \, (\nu_2)_u}{\gamma_2
(\nu_1-\nu_2)}\, l},\\
[4mm]
Y_u  = & - \, \displaystyle{\frac{\gamma_1 \, (\nu_2)_u}
{\gamma_2(\nu_1-\nu_2)}\, X}, & & \\
[4mm]
l_u  = & - \, \displaystyle{\frac{\nu_1 \, (\nu_2)_u}
{\gamma_2(\nu_1-\nu_2)}}\, X;  & &\\
[5mm]
X_v = & &\displaystyle{\frac{\gamma_2 \, (\nu_1)_v}
{\gamma_1(\nu_1-\nu_2)}\, Y,} &\\
[3mm]
Y_v  = &- \, \displaystyle{\frac{\gamma_2 \, (\nu_1)_v}
{\gamma_1(\nu_1-\nu_2)}}\, X & &
+ \, \displaystyle{\frac{\nu_2 \,(\nu_1)_v}{\gamma_1
(\nu_1-\nu_2)}}\, l,\\
[4mm]
l_v  = &  & - \, \displaystyle{\frac{\nu_2 \, (\nu_1)_v}
{\gamma_1(\nu_1-\nu_2)}}\,Y. &
\end{array}
\leqno(3.2)$$
Finding the integrability conditions of the system (3.2) we can formulate
the fundamental Bonnet theorem in the following invariant form:
\begin{thm}\label{T:2.1}
Given four functions $\nu_1(u,v), \, \nu_2(u,v), \, \gamma_1(u,v),
\, \gamma_2(u,v); \; (u,v)\in \mathcal D$ satisfying the following
conditions:
$$\begin{array}{ll}
1) & \nu_1-\nu_2 > 0, \quad \gamma_1 \, (\nu_1)_v > 0, \quad
\gamma_2 \, (\nu_2)_u > 0; \\
[4mm]
2.1) & \displaystyle{\left(\ln \frac{(\nu_1)_v}{\gamma_1} \right)_u=
\frac{(\nu_1)_u}{\nu_1-\nu_2},}
\qquad
\displaystyle{\left(\ln\frac{(\nu_2)_u}{\gamma_2}\right)_v=
-\frac{(\nu_2)_v}{\nu_1-\nu_2};}\\
[5mm]
2.2) & \displaystyle{\frac{\nu_1-\nu_2}{2}\left(\frac{(\gamma_1^2)_v}
{(\nu_1)_v}-\frac{(\gamma_2^2)_u}{(\nu_2)_u}\right)-
(\gamma_1^2+\gamma_2^2)=\nu_1\nu_2.}
\end{array}$$
Then there exists a unique (up to a motion) strongly regular surface
with prescribed invariants
$\nu_1, \, \nu_2, \, \gamma_1, \, \gamma_2$.
\end{thm}
\emph{Proof:} Let $X(u,v), \, Y(u,v), \, l(u,v)$ be three unknown
vector valued functions.
We consider the system (3.2), which can be written in the form
$$\left(
\begin{array}{l}
X_u\\
[2mm]
Y_u\\
[2mm]
l_u \end{array} \right)=
A \left(
\begin{array}{l}
X\\
[2mm]
Y\\
[2mm]
l  \end{array} \right), \qquad
\left(
\begin{array}{l}
X_v\\
[2mm]
Y_v\\
[2mm]
l_v \end{array} \right)=
B \left(
\begin{array}{l}
X\\
[2mm]
Y\\
[2mm]
l  \end{array} \right),\leqno(3.3)
$$
where $A$ and $B$ are given skew symmetric $3\times3$ matrices.

The integrability condition of the system (3.3) is given by the equality
$$B_u-A_v=[A,B]. \leqno(3.4)$$

Applying (3.4) to (3.2), we obtain that (3.4) is equivalent to the
conditions 2.1) and 2.2) of the theorem.

Now, let $x_0X_0Y_0l_0$ be an initial orthonormal right oriented
coordinate system. Applying the theorem of existence and uniqueness
of a solution to (3.3) with initial conditions
$$X(u_0,v_0)=X_0, \quad Y(u_0,v_0)=Y_0, \quad l(u_0,v_0)=l_0,$$
we obtain a unique solution $X(u,v), \, Y(u,v), \, l(u,v); \;
(u,v) \in \mathcal D'$, $(u_0,v_0)\in \mathcal D'\subset \mathcal D$.
\vskip 2mm
Next we prove that $XYl$ form an orthonormal right oriented frame
field in $\mathcal D'$.

Let $X(X^1,X^2,X^3), \, Y(Y^1,Y^2,Y^3), \, l(l^1,l^2,l^3)$ and set
$$f^{ij}(u,v):=X^iX^j+Y^iY^j+l^il^j, \quad i,j=1,2,3.\leqno(3.5)$$
Differentiating (3.5), and taking into account (3.2), we obtain
$f^{ij}(u,v)={\rm const}=\delta_{ij}$, $\delta_{ij}$ being the
Kronecker's deltas. This proves that $XYl$ form an orthonormal
right oriented frame field at any point $(u,v) \in \mathcal D'$.

Further, let us put
$$x_u=\frac{(\nu_2)_u}{\gamma_2(\nu_1-\nu_2)}\,X, \qquad
x_v=\frac{(\nu_1)_v}{\gamma_1(\nu_1-\nu_2)}\,Y.\leqno(3.6)$$
It follows from the conditions of the theorem that the integrability
conditions for (3.6) are valid. Hence, there exists a unique surface
$S$ in a neighborhood $\mathcal D_0$ of $(u_0,v_0)$
$(\mathcal D_0 \subset \mathcal D')$ satisfying the initial condition
$x(u_0,v_0)=x_0$.

The vector valued functions $X, \, Y, \, l$ satisfy (3.2) and in view
of (3.6) the functions $E=x_u^2, \; G=x_v^2$ satisfy (3.1), which shows
that the invariants of $S$ are the given functions
$\nu_1, \, \nu_2, \, \gamma_1, \, \gamma_2$.
\hfill{\qed}
\vskip 2mm
\emph{Geometric illustration of Theorem $3.2$.} Next we give the steps
of a geometric construction of a surface with prescribed invariants.

1. The principal line $c_1: \, x_1=x_1(u,v), \, u\in J_1, \, v=0$ of
the family ${\mathcal F}_1$ with curvature $\kappa_1$, torsion $\tau_1$
and arc-length $s_1$ has the following natural equations:
$$\begin{array}{l}
\kappa_1=\sqrt{\gamma_1^2+\nu_1^2},\\
[2mm]
\tau_1=\displaystyle{-\frac{\gamma_2(\nu_1-\nu_2)}{(\nu_2)_u}
\left(\arctan \frac{\gamma_1}{\nu_1}\right)_u},\\
[4mm]
\displaystyle{\frac{d s_1}{d u}=\sqrt E=\frac{(\nu_2)_u}
{\gamma_2(\nu_1-\nu_2)}}.\end{array}$$

Then there exists a unique (up to a motion) curve
$c_1: \, x_1=x_1(u), \, u\in J_1$ with prescribed $\kappa_1, \tau_1$
and arc-length $s_1$.

2. There exists a unique orthonormal frame field $X(u), \, Y(u), \, l(u)$
along $c_1$ with torse forming normals $\{Y, l\}$ satisfying the
initial conditions
$$X(0)=X_0, \; Y(0)=Y_0, \; l(0)=l_0.$$

3. Any principal curve $c_2: \; x_2=x(u,v), \, u={\rm const}$ of
the family ${\mathcal F}_2$ with curvature $\kappa_2$, torsion
$\tau_2$ and arc-length $s_2$ has the following natural equations:
$$\begin{array}{l}
\kappa_2=\sqrt{\gamma_2^2+\nu_2^2},\\
[2mm]
\tau_2=\displaystyle{-\frac{\gamma_1(\nu_1-\nu_2)}{(\nu_1)_v}
\left(\arctan \frac{\gamma_2}{\nu_2}\right)_v},\\
[4mm]
\displaystyle{\frac{d s_2}{d v}=\sqrt G=\frac{(\nu_1)_v}
{\gamma_1(\nu_1-\nu_2)}}.\end{array}$$

These equations completely determine the curve $c_2$ with the initial
orthonormal coordinate system $x_1(u)X(u)Y(u)l(u)$.

Hence, the surface is determined uniquely (up to a motion) by the
four invariants $\nu_1, \, \nu_2; \, \gamma_1,$ \, $\gamma_2$.

\section{Weingarten surfaces and geometric parameters}
We shall use the definition for the class of Weingarten surfaces
in the following form  \cite{W1,W2}.

A surface $S$ is said to be a \emph{Weingarten surface} if there exists
a differentiable function $w(t), \; t\in \mathcal{I}\subseteq{\R}$
such that the principal curvatures of $S$ at every point satisfy the
condition $\nu_2=w(\nu_1)$.

In section 2 we described the class of surfaces satisfying the condition
$\gamma_1=0$. Now we shall specify this result for the class of
Weingarten surfaces.

First we recall some basic facts about canal surfaces.

A surface $S$ is said to be \emph{canal} if it is an envelope of a
one-parameter family of spheres.

Any canal surface $S$ is a one-parameter family of circles
$S=\{c(u)\}_{u\in J}$ and has the following properties \cite {K}, (p. 33):
\vskip 2mm
\emph{$S =\{c(u)\}_{u\in J}$ is a canal surface if and only if the
tangent planes $($normals$)$ to $S$ at the points of the circles
$c(u), \, u={\rm const}$ pass through a fixed point
$z(u)\; (\tilde z(u))$.}
\vskip 2mm
\emph{If $S =\{c(u)\}_{u\in J}$ is a canal surface, then the one family of principal
lines of $S$ consists of the circles $\{c(u)\}_{u\in J}$.}
\vskip 2mm
The inverse statement of the last assertion is also true. It has been proved in \cite{G}
(Theorem 20.12) under the assumption that the focal set of $S$ is a
differentiable manifold. We need this property in the following form.

\begin{prop}\label {P:4.1}
If one of the families of principal lines of a surface $S$ consists of
circles, then $S$ is a canal surface.
\end{prop}

\emph{Proof:} Let the family ${\mathcal F}_2$ consist of circles, i.e.
$\tau_2=0$ and $\kappa_2 =\kappa_2(u)={\rm const}$ ($u$ - fixed).
Then (2.8) implies that $\theta_2=\theta_2(u)={\rm const}$ ($u$ - fixed).
Hence $S$ is a canal surface.
\hfill{\qed}
\vskip 2mm
Let $S: \; x=x(u,v), \;(u,v)\in \mathcal D$ be a canal surface, i.e.
the envelope of a family of spheres $\{S^2(z(u), R(u))\}, \, u\in J$,
where $z(u)$ and $R(u)$ are the center and the radius of the sphere
$S^2(u), \; u\in J$, respectively. We assume that the family
${\mathcal F}_2$ of principal lines consists of the intersections
$S^2(u)\cap S, \; u\in J$, which are circles.

As usual,
$Y$ is the unit tangent vector field, normal to the circles of
$\mathcal F_2$. We denote by $g$ the induced metric on $S$ and by
$\eta$ the 1-form associated with $Y$, i.e. $\eta(Z)=g(Y,Z)$ for any
tangent vector field $Z$ to $S$. Adapting the formulas from \cite{GM}
to the case $n=2$, we obtain that the second fundamental form of any
canal surface $S$ has the following form:
$$h=\frac{1}{R}\,g-\frac{1-{R'}^2}{R(1-{R'}^2-RR''-R\sqrt{1-{R'}^2}
\,k \cos v)}\,\eta\otimes \eta,$$
where $k$ is the curvature of the curve of the centers
$z=z(u),\;u\in J$ and $v\in [0, 2\pi)$.

Then the principal curvatures of $S$ are given as follows:
$$\nu_1=\frac{1}{R}, \quad \nu_2=-\frac{R''+\sqrt{1-{R'}^2}\, k \cos v}
{1-{R'}^2-RR''-R\sqrt{1-{R'}^2}\, k \cos v}.\leqno(4.1)$$

Now we can prove the following modification of Theorem 2.2 for the case
of Weingarten surfaces.
\begin{thm}\label {T:4.1}
Any Weingarten surface satisfying the condition
$\gamma_1=0$ is locally part of a rotational surface.
\end{thm}

\emph{Proof:} Let $S$ be parameterized by principal parameters and assume that
the family $\mathcal F_1$ consists of geodesics, i.e. $\gamma_1=0$.
Then the formulas (2.3), (2.5) and (2.6) immediately imply that
$\nu_1=\nu_1(u),\; \kappa_1=\kappa_1(u)$ and $\tau_1= 0$.

If $S$ is a Weingarten surface, then $\nu_2=w(\nu_1)$ and
$\nu_2=\nu_2(u)$. Taking into account (2.3) and $\nu_1-w(\nu_1)\neq 0$, we find that
$$\gamma_2=\frac{w'\,X(\nu_1)}{\nu_1-w(\nu_1)}$$
and therefore $\gamma_2=\gamma_2(u)$. Because of (2.7) and (2.8) it
follows that  $\tau_2= 0$ and $\kappa_2=\kappa_2(u).$

Applying Proposition 4.1 we obtain that $S$ is a canal surface. Then
(4.1) in view of the condition $\nu_2=\nu_2(u)$ implies that
$k=0$, i.e. the curve of centers $z=z(u)$ is part of a straight line. Hence
$S$ is a rotational surface.
\hfill{\qed}
\vskip 2mm
\begin{rem}
The argument in the above statement implies that in terms of the
denotations in Section 2 the class of rotational surfaces is characterized
by the conditions:
$$\gamma_1=0, \quad \kappa ={\rm const}, \quad \tau=0.$$
\end{rem}

>From now on we consider only strongly regular Weingarten surfaces. For the sake of
symmetry with respect to the principal curvatures $\nu_1$ and $\nu_2$ we shall use
the following characterization of strongly regular Weingarten surfaces:

\emph{A strongly regular surface $S:\; x=x(u,v),\; (u,v)\in \mathcal{D}$ is Weingarten
if there exist two differentiable functions
$f(\nu), \; g(\nu), \; f(\nu)-g(\nu) > 0, \; f'(\nu)g'(\nu)\neq 0, \;
\nu \in \mathcal{I}\subseteq{\R}$
such that the principal curvatures of $S$ at every point are given by
$\nu_1=f(\nu), \; \nu_2=g(\nu), \; \nu=\nu(u,v), \; \nu_u(u,v)\nu_v(u,v) \neq 0, \;
(u,v) \in \mathcal D$.}
\begin{lem}\label {L:4.1} Let $S:\; x=x(u,v),\; (u,v)\in \mathcal{D}$ be a strongly
regular Weingarten surface parameterized with principal parameters. Then the function
$$\lambda = \ln\left\{\sqrt E \exp\left(\int \frac{f'd\nu}{f-g}\right)\right\}$$
does not depend on $v$, while the function
$$\mu = \ln \left\{\sqrt G \exp \left(\int \frac{g'd\nu}{g-f}\right)\right\}$$
does not depend on $u$.
\end{lem}

{\it Proof}: Taking into account (2.3) and (2.1), we find
$$\gamma_1=\frac{f'(\nu) Y(\nu)}{f(\nu)-g(\nu)}=-Y(\ln \sqrt E), \qquad
\gamma_2=\frac{g'(\nu) \,X(\nu)}{f(\nu)-g(\nu)}=
X(\ln \sqrt G),$$
which imply that
$$Y\left(\int\frac{f'(\nu)\,d \nu}{f(\nu)-g(\nu))} + \ln \sqrt E\right) = 0,
\qquad
X\left(\int\frac{g'(\nu)\, d \nu}{g(\nu)-f(\nu)} + \ln \sqrt G\right) = 0.$$
The last equalities mean that $\lambda_v=0$ and $\mu_u=0$. \qed
\vskip 2mm
\begin{rem}
In the case of a surface with Gaussian curvature $K=-1$ in \cite{G} (Lemma 21.8) it is
proved that $\lambda=\lambda(u)$ and $\mu=\mu(v)$.
\end{rem}
\vskip 2mm
Now we introduce the notion of geometric principal parameters on a strongly rergular
Weingarten surface.

\begin{defn}
Let $S:\; x=x(u,v),\; (u,v)\in \mathcal{D}$ be a strongly regular Weingarten surface
parameterized with principal parameters. The parameters $(u, v)$ are said to be
\emph{geometric} principal, if the functions $\lambda(u)$ and $\mu(v)$ from Lemma \ref{L:4.1}
are constants.
\end{defn}
\begin{thm}\label {T:4.2}
Any strongly regular Weingarten surface admits locally geometric principal
parameters.
\end{thm}
\emph{Proof}: Let $S:\; x=x(u,v),\; (u,v)\in \mathcal{D}$ be a strongly
regular Weingarten surface parameterized with principal parameters. Then
$\nu_1=f(\nu), \; \nu_2=g(\nu), \; \nu=\nu(u,v)$ for some differentiable functions
$f$, $g$ and $\nu$ satisfying the conditions
$$f(\nu)-g(\nu)> 0,\quad f'(\nu)g'(\nu)\neq 0, \quad \nu_u(u,v)\nu_v(u,v) \neq 0;
\quad (u,v) \in \mathcal D.$$

Let $a={\rm const} \neq 0, \; b={\rm const} \neq 0$, \;
$(u_0, v_0) \in \mathcal D$ and $\nu_0=\nu(u_0,v_0)$.
We change the parameters $(u, v)\in \mathcal{D}$ with
$(\bar u, \bar v)\in \bar{\mathcal{D}}$ by the formulas
$$\begin{array}{l}
\displaystyle{\bar u=a\int_{u_0}^u \sqrt E \exp\left(\int_{\nu_0}^{\nu} \frac{f'd\nu}{f-g}\right)}
\, du\, +\overline{u}_0 , \quad \bar u_0 = {\rm const},\\
[4mm]
\displaystyle{\bar v=b\int_{v_0}^v\sqrt G \exp \left(\int_{\nu_0}^{\nu} \frac{g'd\nu}{g-f}\right)\,
dv}\,+\overline{v}_0, \quad \bar v_0={\rm const}. \end{array}$$
According to Lemma \ref{L:4.1} it follows that $(\bar u, \bar v)$ are again principal parameters and
$$\bar E=\frac{1}{a^2}\,\exp \left(-2\int_{\nu_0}^{\nu} \frac{f'd\nu}{f-g}\right),\quad
\bar G=\frac{1}{b^2}\,\exp \left(-2\int_{\nu_0}^{\nu} \frac{g'd\nu}{g-f}\right).\leqno(4.2)$$
For the functions from Lemma \ref{L:4.1} we find
$$\bar \lambda =1/a^2, \quad \bar \mu =1/b^2.$$
Furthermore $1/a^2=\bar E(u_0,v_0), \; 1/b^2 = \bar G(u_0,v_0)$. {\qed}
\vskip 2mm
>From now on we assume that the considered strongly regular Weingarten surface
$S: x(u,v), \; (u,v) \in \mathcal D$ is parameterized with geometric principal parameters $(u, v)$.

As an immediate consequence from Theorem \ref{T:4.2}, we get
\begin{cor}
Let $S$ be a strongly regular Weingarten surface parameterized by geometric principal
parameters $(u, v)$. Then any geometric principal parameters $(\tilde u, \tilde v)$
on $S$ are determined by $(u, v)$ up to an affine transformation of the type
$$\begin{array}{lccc}
\tilde u = a_{11} \, u  & & +\, b_1,\\
[2mm]
\tilde v = & a_{22} \, v &+\, b_2, \end{array}
\quad a_{11}a_{22} \neq 0,
$$
or
$$\begin{array}{lcc}
\tilde u = &a_{12} \, v & + \,c_1,\\
[2mm]
\tilde v = a_{21} \, u & &+\,c_2, \end{array}
\quad a_{12}a_{21}\neq 0, $$
where $a_{ij}, \, b_i, \, c_i; \, i,j=1,2 $ are constants.
\end{cor}
Next we give a simple criterion principal parameters to be geometric.
\begin{prop}\label{P:4.2}
Let a strongly regular Weingarten surface $S:\; x=x(u,v),\; (u,v)\in \mathcal{D}$ be
parameterized with principal parameters. Then $(u,v)$ are geometric if and only if
$$\sqrt{EG}(\nu_1-\nu_2)={\rm const}.$$
\end{prop}

{\it Proof}: The equality $\sqrt{EG}\,|\nu_1-\nu_2|=\exp(\lambda + \mu)$ and
Lemma \ref{L:4.1} imply the assertion. \qed
\vskip 2mm
Taking into account Theorem \ref {T:4.2} and (3.1), we obtain the following statement.
\begin{cor}\label {C:4.1}
A strongly regular Weingarten surface admits isothermal geometric
principal parameters if and only if
$f'(\nu)+g'(\nu)=0$, i. e. $H=\displaystyle{\frac{f(\nu)+g(\nu)}{2}}={\rm const}$.
\end{cor}

\begin{thm}\label {T:4.3}
Given two differentiable functions $f(\nu), \, g(\nu); \; \nu \in \mathcal{I},$ \;
$f(\nu)-g(\nu)> 0$, \; $f'(\nu)g'(\nu)\neq 0$ and a differentiable function
$\nu(u,v), \; (u,v) \in {\mathcal D}$ satisfying the conditions
$$\nu_u\,\nu_v\neq 0, \quad \nu(u,v)\in \mathcal{I}.$$

Let $(u_0, v_0) \in \mathcal D, \; \nu_0=\nu(u_0, v_0)$ and $a\neq 0, \, b\neq 0$ be two constants.
If
$$\begin{array}{l}
\ds{\;\;\;b^2\exp\left(2\int_{\nu_0}^{\nu}
\frac{g'd\nu}{g-f}\right)\left[f'\nu_{vv}+
\left(f''-\frac{2f'^2}{f-g}\right)\nu^2_v\right]}\\
[4mm] \ds{-a^2\exp\left(2\int_{\nu_0}^{\nu}
\frac{f'd\nu}{f-g}\right) \left[g'\nu_{uu} +
\left(g''-\frac{2g'^2}{g-f}\right)\nu^2_u\right]=fg(f-g)},
\end{array}\leqno(4.3)$$
then there exists a unique (up to a motion) strongly regular Weingarten surface
$S:\; x=x(u,v),$ \, $(u,v)\in \mathcal D_0 \subset \mathcal D$ with invariants
$$\nu_1=f(\nu), \quad \nu_2=g(\nu), $$
$$\gamma_1=\exp\left(\int_{\nu_0}^{\nu} \frac{g'd\nu}{g-f}\right)\,\frac{bf'}{f-g}\,\nu_v,
\quad \gamma_2=-\exp\left(\int_{\nu_0}^{\nu} \frac{f'd\nu}{f-g}\right)\,\frac{ag'}{g-f}\,\nu_u ,
\leqno(4.4) $$
and furthermore $(u,v)$ are geometric principal parameters for $S$.
\end{thm}

\emph{Proof:} Applying Theorem \ref{T:2.1} and taking into account Theorem
\ref {T:4.2}, we obtain that the integrability conditions 2.1) and 2.2)
in Theorem \ref{T:2.1} reduce to (4.3), which proves the assertion.
\hfill{\qed}
\vskip 2mm
Using (2.6) and (4.4), we also get
$$\begin{array}{l}
\ds{\tau_1=\frac{-bf}{f^2+\gamma_1^2}\left\{f'\nu_{uv}+\left[f''- \frac{f'^2(2f-g)}{f(f-g)}
\right]\nu_u\,\nu_v\right\}},\\
[4mm]
\ds{\tau_2=\frac{ag}{g^2+\gamma_2^2}\left\{g'\nu_{uv}+\left[g''- \frac{g'^2(2g-f)}{g(g-f)}
\right]\nu_u\,\nu_v\right\}}\end{array} \leqno(4.5)$$
with respect to the geometric principal parameters $(u, v)$.

\section{Applications to special Weingarten surfaces}

\subsection{Surfaces of constant mean curvature}

Let $S$ be a strongly regular surface of constant mean curvature $H$, whose parametric net is principal.

Putting $2\nu=\nu_1-\nu_2 > 0 $, we have
$$\nu_1=H+\nu, \qquad \nu_2=H-\nu, \qquad K=H^2-\nu^2.$$

Applying Theorem \ref{T:4.2} we can assume that $S$ is parameterized with geometric principal
parameters $(u,v)$. Taking into account (4.2) and Corollary \ref{C:4.1}, we find
$$E=\frac{\nu_0}{a^2\nu}, \quad G=\frac{\nu_0}{a^2\nu}, \quad
L=\frac{\nu_0(H+\nu)}{a^2\nu}, \quad N=\frac{\nu_0(H-\nu)}{a^2\nu}. $$

Choosing $a^2= 2\nu_0$, we obtain special geometric principal parameters $(u, v)$, which we call
\emph{canonical} principal. With respect to these canonical principal parameters we have
$$E=G=\frac{1}{2\nu}, \quad \gamma_1=(\sqrt{2\nu}\,)_v\neq 0, \quad \gamma_2=-(\sqrt{2\nu}\,)_u\neq 0$$
and the Gauss equation (4.3) in Theorem 4.8 gets the following simple form
$$\nu_{uu}+\nu_{vv}=\frac{1}{\nu}\,(\nu_u^2+\nu_v^2)+
(H^2-\nu^2),\leqno(5.1)$$
which can be written as follows:
$$\Delta \ln \nu=\frac{1}{\nu}\,(H^2-\nu^2).\leqno(5.2)$$

Applying Theorem 4.8 we obtain
\begin{cor}\label{C:5.1}
Given a constant $H$ and a function $\nu (u,v) > 0$ in a neighborhood
$\mathcal D$ of $(u_0, v_0)$ with $\nu_u \nu_v \neq 0$,
satisfying the equation $(5.2)$ and an initial orthonormal frame $x_0X_0Y_0l_0$.

Then there exists a unique strongly regular surface
$S: \; x=x(u, v), \; (u, v) \in \mathcal D_0 \; ((u_0, v_0) \in \mathcal D_0
\subset \mathcal D)$, such that

(i) $(u, v)$ are canonical principal parameters;

(ii) $x(u_0, v_0)=x_0, \; X(u_0, v_0)=X_0, \; Y(u_0, v_0)=Y_0, \; l(u_0, v_0)=l_0$;

(iii) $S$ is of constant mean curvature $H$ with invariants

$$\nu_1=H+ \nu, \quad \nu_2=H- \nu,
\quad \gamma_1=(\sqrt {2\nu}\,)_v,  \quad
\gamma_2= -(\sqrt{2\nu}\,)_u.$$
\end{cor}

Up to a renumbering of the principal parameters we can assume that $H \geq 0$.

The equation (5.2) can be specified as follows:
\vskip 2mm
\centerline{$H=0$ \textbf{(Minimal surfaces)}}
\vskip 2mm
In this case (5.2) becomes
$$\Delta \ln \nu = -\nu. \leqno (5.3)$$
Putting $\nu =e^{\lambda}$, the equation $(5.3)$ gets the form
$$\Delta \lambda =-e^{\lambda}. \leqno(5.4)$$

\begin{cor}
Any solution $\lambda$ of the equation $(5.4)$, satisfying the condition
$\lambda_u \lambda_v\neq 0$, generates a unique (up to a motion) minimal strongly
regular surface with invariants
$$\nu_1= e^{\lambda}, \quad \nu_2=-e^{\lambda},
\quad \gamma_1=\frac{1}{\sqrt 2}\, e^{\frac{\lambda}{2}}\,\lambda_v, \quad
\gamma_2= -\frac{1}{\sqrt 2}\,e^{\frac{\lambda}{2}}\,\lambda_u.$$
\end{cor}
\vskip 2mm
\centerline{$H={\rm const}\, >0$ \textbf{(CMC-surfaces)}}
\vskip 2mm
In this case the equation (5.2) can be written in the form
$$\Delta \ln \nu=-2H\,\frac{\left(\frac{\nu}{H}\right)-
\left(\frac{H}{\nu}\right)}{2}.\leqno(5.5)$$
Putting $\displaystyle{\frac{\nu}{H}= e^{\lambda}}$,
then $(5.5)$ becomes
$$\Delta \lambda =-2 H\,\sinh \lambda. \leqno(5.6)$$
Especially, if $H=\ds{\frac{1}{2}}$, then
$$\Delta \lambda = - \sinh \lambda.
\leqno(5.7)$$

Corollary \ref{C:5.1} implies the following statement \cite{KO,W} (See also \cite{KK}):

\begin{cor}
Any solution $\lambda$ of the equation $(5.7)$, satisfying the condition
$\lambda_u \lambda_v\neq 0$, generates a unique (up to a motion) strongly regular
surface of constant mean curvature $1/2$ with invariants
$$\nu_1=\ds{\frac{1}{2}(1 +  e^{\lambda}), \quad
\nu_2=\frac{1}{2}(1 -  e^{\lambda})},
\quad \gamma_1=\frac{1}{\sqrt{2}}\, e^{\frac{\lambda}{2}}\,\lambda_v, \quad
\gamma_2= -\frac{1}{\sqrt{2}}\,e^{\frac{\lambda}{2}}\,\lambda_u.$$
\end{cor}

\subsection{Surfaces of constant Gaussian curvature}

Let $S$ be a strongly regular surface, whose parametric net is principal. In this subsection
we consider surfaces with Gaussian curvature $K=\pm 1$.

Putting $\nu=\nu_1$, we have
$$\nu_2=\frac{K}{\nu}, \quad \nu_1-\nu_2= \frac{\nu^2-K}{\nu}.$$

Further, we assume that the surface $S$ is parameterized by geometric principal
parameters and $\nu_1-\nu_2 > 0$. Then
$$E=\frac{\nu_0^2-K}{a^2(\nu^2-K)}, \quad G= \frac{\nu_0^2-K}{b^2\nu_0^2}
\frac{\nu^2}{\nu^2-K}.$$

We denote $\varepsilon = {\rm sign} \, \nu={\rm sign \, (\nu^2-K)}$. Choosing
$a^2=\varepsilon (\nu_0^2-K)$ and $b^2=\ds{\varepsilon \frac{\nu_0^2-K}{\nu_0^2}}$,
we obtain special geometric principal parameters
$(u, v)$, which we call \emph{canonical} principal. With respect to these canonical
principal parameters we have
$$ E=\frac {\varepsilon}{\nu^2-K}, \quad G=\frac{\varepsilon \nu^2}{\nu^2-K}; \quad
\gamma_1=\frac{\nu_v}{\sqrt{\varepsilon(\nu^2-K)}}, \quad
\gamma_2=\frac{-K\nu_u}{\nu\sqrt{\varepsilon(\nu^2-K)}}.$$
The Gauss equation (4.5) in Theorem 4.8 gets the form
$$\nu_{vv}+K\nu_{uu}-\frac{2\nu}{\nu^2-K}(\nu_v^2+K\nu_u^2)=K\nu.\leqno(5.8)$$

Applying Theorem 4.8 we obtain
\begin{cor}
Given a constant $K = \pm 1$ and a function $\nu(u,v)$ in a
neighborhood $\mathcal D$ of $(u_0, v_0)$ with $\nu (\nu^2-K) > 0, \;
\nu_u \nu_v \neq 0$, satisfying the
equation $(5.8)$, and an initial orthonormal frame $x_0X_0Y_0l_0$.

Then there exists a unique strongly regular surface $S: \; x=x(u, v), \; (u, v) \in \mathcal D_0
\; ((u_0, v_0) \in \mathcal D_0 \subset \mathcal D)$,
such that

(i) $(u, v)$ are canonical principal parameters;

(ii) $x(u_0, v_0)=x_0, \; X(u_0, v_0)=X_0, \; Y(u_0, v_0)=Y_0, \; l(u_0, v_0)=l_0$;

(iii) $S$ is of constant Gaussian curvature $K$ with invariants

$$\nu_1=\nu, \quad \nu_2=\frac{K}{\nu}, \quad \gamma_1=
\frac{\nu_v}{\sqrt{\varepsilon(\nu^2-K)}}, \quad
\gamma_2=\frac{-K\nu_u}{\nu\sqrt{\varepsilon(\nu^2-K)}},\quad
\varepsilon = {\rm sign \, (\nu^2-K)}.$$
\end{cor}

The equation (5.8) can be specified as follows:
\vskip 2mm
\centerline{\bf The case $K=1$}
\vskip 2mm
Let $\varepsilon = 1$. Then the equation (5.8) becomes
$$\Delta \ln \frac{\nu-1}{\nu+1}=\frac{2\nu}{\nu^2-1},\leqno(5.9)$$
or putting $\displaystyle{\nu=-\coth \frac{\lambda}{2}}$,
$$\Delta \lambda =-\sinh \lambda, \quad \lambda<0. \leqno(5.10)$$

\begin{cor}
Any solution $\lambda$ of the equation $(5.10)$, satisfying the conditions
$\lambda < 0, \; \lambda_u \lambda_v\neq 0$ generates a unique (up to a motion) strongly
regular surface with Gaussian curvature $K=1$ and invariants
$$\nu_1=-\coth \frac{\lambda}{2}, \quad \nu_2=-\tanh \frac{\lambda}{2}, \quad
\gamma_1=-\frac{\lambda_v}{2\sinh \frac{\lambda}{2}}\, , \quad
\gamma_2=-\frac{\lambda_u}{2\cosh \frac{\lambda}{2}}\,.$$
\end{cor}
Let $\varepsilon = -1$. Renumbering the parameters, we get the previous case.
\vskip 2mm
\centerline{\bf The case $K=-1$}
\vskip2mm
In this case $\varepsilon = 1$ and the equation (5.8) becomes
$$\Delta^h(\arctan \nu)=\frac{\nu}{1+\nu^2},\leqno(5.11)$$
where $\Delta^h$ denotes the hyperbolic Laplacian.
Putting $\nu=\tan {\frac{\lambda}{2}}$,
the equation (5.11) gets the form
$$\Delta^h \lambda = \sin \lambda,
\quad 0<\lambda<\pi.\leqno(5.12)$$

\begin{cor}
Any solution $\lambda$ of the equation $(5.12)$, satisfying the conditions
$0< \lambda < \displaystyle{\frac{\pi}{2}}, \; \lambda_u\lambda_v\neq 0$,
generates a unique (up to a motion) strongly regular surface with Gaussian curvature $K=-1$
and invariants
$$\nu_1=\tan{\frac{\lambda}{2}}, \quad \nu_2=-\cot{\frac{\lambda}{2}},
\quad \gamma_1=\frac{\lambda_v}{2\cos{\frac{\lambda}{2}}}, \quad
\gamma_2=\frac{\lambda_u}{2\sin{\frac{\lambda}{2}}}.$$
\end{cor}
\vskip 2mm
\begin{rem}
Let a strongly regular surface with Gaussian curvature $K=-1$ be parameterized with
canonical principal parameters. Then $E+G=1$. According to Lemma 21.5 in \cite{G}
it follows that the parametric net is a Tchebyshef net.
\end{rem}
{\bf Example. The Kuen surface with respect to canonical principal parameters.}
\vskip 3mm

The Kuen surface with respect to principal parameters is given as follows \cite{G}:
$$S: \quad
\left\{\begin{array}{l}
\displaystyle{x_1=2(\cos\,u+u\,\sin\,u)\frac{\sin\,v}{1+u^2\sin^2 v},}\\
[4mm]
\displaystyle{x_2=2(\sin\,u-u\,\cos\,u)\frac{\sin\,v}{1+u^2\sin^2 v},}\\
[4mm]
\displaystyle{x_3=\ln{\tan{\left(\frac{v}{2}\right)}} + 2\frac{\cos\,v}{1+u^2\sin^2 v};}
\end{array}\right. \qquad u\in (0, 2\pi), \; v\in \left(0, \frac{\pi}{2}\right).
$$
The coefficients of the first and the second fundamental form are the following:
$$E=\frac{4u^2\sin^2 v}{(1+u^2\sin^2 v)^2}, \quad F=0, \quad G=\frac{(1-u^2\sin^2 v)^2}{\sin^2v
(1+u^2\sin^2v)^2};$$

$$L=-2u\sin\,v \frac{(1-u^2\sin^2 v)}{(1+u^2\sin^2v)^2}, \quad M=0, \quad
N=\frac{2u}{\sin\,v}\frac{(1-u^2\sin^2 v)}{(1+u^2\sin^2v)^2}.$$

The invariants of $S$ are
$$\nu_1=-\frac{1-u^2\sin^2v}{2u\,\sin\,v}, \qquad \nu_2=\frac{2u\,\sin\,v}{1-u^2\sin^2v};$$

$$\gamma_1=-\cos \,v, \qquad \gamma_2=\frac{-2\,\sin\,v}{1-u^2\sin^2v}.$$
\vskip 2mm
We parameterize the Kuen surface with canonical principal parameters:
$$\bar u=u,\quad \bar v=\ln{\tan{\left(\frac{v}{2}\right)}}.\leqno(5.13)$$

The Kuen surface with respect to canonical principal parameters (again denoted by $(u,v)$)
has the following representation:
$$S: \quad
\left\{ \begin{array}{l}
\displaystyle{x_1=2(\cos\,u+u\,\sin\,u)\frac{\cosh\,v}{u^2+\cosh^2 v},}\\
[4mm]
\displaystyle{x_2=2(\sin\,u-u\,\cos\,u)\frac{\cosh\,v}{u^2+\cosh^2 v},}\\
[4mm]
\displaystyle{x_3=v - \frac{\sinh{2v}}{u^2+\cosh^2 v}.}
\end{array}\right.$$

The corresponding coefficients and invariants are as follows
$$E=\frac{4u^2\cosh^2 v}{(u^2+\cosh^2 v)^2}, \quad F=0, \quad G=\frac{(u^2-\cosh^2 v)^2}
{(u^2+\cosh^2v)^2};$$
$$L=2u\cosh\,v \frac{(u^2-\cosh^2 v)}{(u^2+\cosh^2v)^2}, \quad M=0, \quad
N=-2u\,\cosh\,v \frac{(u^2-\cosh^2 v)}{(u^2+\cosh^2v)^2};$$

$$\nu_1=\frac{u^2-\cosh^2v}{2u\,\cosh\,v},\qquad \nu_2=-\frac{2u\,\cosh\,v}{u^2-\cosh^2v};$$

$$\gamma_1=\tanh \,v,\qquad \gamma_2=\frac{2\,\cosh\,v}{u^2-\cosh^2v}.$$

\begin{rem} (Geometric interpretation of the canonical principal parameters on the Kuen surface)
Taking into account Theorem 21.30 in \cite{G}, it follows that the canonical principal parameters
(5.13) are the Bianchi transform of the usual principal parameters of the pseudo-sphere.
Furthermore this parametric net is a Tchebyshef net.
\end{rem}

\end{document}